\documentclass[final]{elsart}
\usepackage{epsfig}
\usepackage{amssymb}
\usepackage{amsmath}

  \def\phi{\varphi}    
\def\ell{{l}}     \def\proj{\text{proj}}
\def\tld{\widetilde}

 \def\RR{{\mathbb{R}}} \def\ZZ{{\mathbb{Z}}} \def\FF{F} 
\def\QQ{{\mathbb{Q}}}

\begin{document}

\begin{frontmatter}
\runauthor{Hanusa}
\title{Pseudo-centrosymmetric matrices, with applications to counting perfect matchings}

\author[C]{Christopher R. H. Hanusa}
\address[C]{Department of Mathematical Sciences,
Binghamton University, Binghamton, New York, 13902-6000 
{\tt chanusa@math.binghamton.edu}}

\begin{abstract}
We consider square matrices $A$ that commute with a fixed square matrix $K$, both with entries
in a field $\FF$ not of characteristic 2.  When $K^2=I$, Tao
and Yasuda defined $A$ to be generalized centrosymmetric with respect to $K$.  When $K^2=-I$, we define $A$ to be
pseudo-centrosymmetric with respect to $K$; we show that the determinant of every even-order pseudo-centrosymmetric
matrix is the sum of two squares over $\FF$, as long as $-1$ is not a square in $\FF$.  When a
pseudo-centrosymmetric matrix $A$ contains only integral entries and is pseudo-centrosymmetric with respect to a
matrix with rational entries, the determinant of $A$ is the sum of two integral squares.
This result, when specialized to when $K$ is the even-order alternating exchange matrix, applies to enumerative
combinatorics.  Using solely matrix-based methods, we reprove a weak form of Jockusch's theorem for enumerating
perfect matchings of 2-even symmetric graphs.  As a corollary, we reprove that the number of domino tilings of
regions known as Aztec diamonds and Aztec pillows is a sum of two integral squares.
\end{abstract}


\end{frontmatter}

\maketitle


\section{Introduction}
\label{Introduction}

This article is divided into two halves; results concerning certain types of matrices from the first half apply to
enumerative combinatorics in the second half.  An outline of the structure of the article follows.

Centrosymmetric matrices have been studied in great detail; in Section \ref{sec:PseudoCentro}, we use their
definition as a starting point to define an extension we call {\em pseudo-centrosymmetric} matrices, defined over a
field $\FF$, not of characteristic 2.  Theorem \ref{SumofSquares} proves that a subclass of these matrices have
determinants that are the sum of two squares in $\FF$. Specializing further in Section \ref{sec:AltCentro}, we show
that when the matrix is {\em alternating centrosymmetric}, its determinant has a nice symmetric form.  Returning to
the general case in Section \ref{sec:general}, Theorem \ref{thm:SumofSquaresGen} establishes that when $-1$ is not
a square in $\FF$, even-order pseudo-centrosymmetric matrices over $\FF$ have determinants that
are a sum of two squares in $\FF$.  A direct corollary is that when the entries are integers, the determinant is a
sum of integral squares.

In Section \ref{sec:Enum}, we apply these results on matrices to the question of counting perfect matchings 
of graphs.  Jockusch proved that the number of perfect matchings of a 2-even symmetric graph is always the 
sum of squares.  In Section \ref{sec:JockuschThm}, we are able to reprove a weaker version of Jockusch's 
theorem using only matrix methods in the case when the graph in question is embedded in the square lattice 
and the center of rotational symmetry is in the center of a unit lattice square.  (A restriction of this 
type is necessary, as shown in Section \ref{sec:CtrEx}.)  Graphs of this restricted type occur in the study 
of domino tilings of nice regions called Aztec diamonds, as explained in Section \ref{sec:Pillows}.

Future directions of study appear in various remarks throughout the paper.

\section{Matrix-Theoretical Results}
\label{sec:Mtx}

Throughout Sections \ref{sec:PseudoCentro} and \ref{sec:AltCentro}, we define our matrices over an arbitrary field
$\FF$, not of characteristic 2.

\subsection{Pseudo-centrosymmetric matrices with respect to $K$}
\label{sec:PseudoCentro}

Define $J$ to be the $n\times n$ {\em exchange matrix} with 1's along the cross-diagonal ($j_{i,n-i+1}$) and 0's
everywhere else.  Matrices $A$ such that $JA=AJ$ are called {\em centrosymmetric} and have been studied in much
detail because of their applications in wavelets, partial differential equations, and other areas (see
\cite{Muir,Weaver}).  A matrix such that $JA=-AJ$ is called {\em skew-centrosymmetric}. Tao and Yasuda define a
generalization of these matrices for any involutory matrix $K$ ($K^2=I$).  A matrix $A$ that is {\em
centrosymmetric with respect to $K$} satisfies $KA=AK$ (see \cite{Andrew, TaoYasuda}).  A matrix $A$ that is
{\em skew-centrosymmetric with respect to $K$} satisfies $KA=-AK$.

In the study of generalized Aztec pillows, a related type of matrix arises.  Define a matrix $K$ to be {\em 
anti-involutory} if $K^2=-I$.  We will call a matrix {\em pseudo-centrosymmetric with respect to $K$} if 
$KA=AK$, or {\em pseudo-skew-centrosymmetric with respect to $K$} if $KA=-AK$.

{\rem Another definition of centrosymmetric matrices is that $KAK=A$.  We must be careful with their 
pseudo-analogues because if $KA=AK$ when $K^2=-I$, then $KAK=-A$.
}

In this article, we focus on studying even-order pseudo- and pseudo-skew-centrosymmetric matrices over 
$\FF$. When $K$ is an 
$n\times n$ matrix for $n=2k$ even, we can write
\begin{equation*}
K=
\left(
\begin{matrix}
K_1 & K_2 \\
K_3 & K_4 \\
\end{matrix}
\right),
\end{equation*}
 with each submatrix $K_i$ being of size $k\times k$.  We will explore more general anti-involutory matrices $K$ in
 Section \ref{sec:general}, but for now we will focus on the simple case when 
$K_1=K_4=0$.  Since $K$ is anti-involutory, $K_3=-K_2^{-1}$, so $K$ is of the form
\begin{equation}
K=
\left(
\begin{matrix}
 0 & K_2 \\
-K_2^{-1} & 0 \\
\end{matrix}
\right).
\label{SimpleK}
\end{equation}
For such a matrix $K$, a matrix $A$ that is a pseudo-centrosymmetric with respect to $K$ has a simple form for 
its determinant.

{\thm
Let $K$ and $A$ be matrices defined over $\FF$.  If $K$ is an anti-involutory matrix of 
the form in Equation \eqref{SimpleK}, and the $2k\times 2k$ matrix $A$ is pseudo-centrosymmetric with 
respect to $K$, then $A$ 
has the form 
\begin{equation*}
A=
\left(
\begin{matrix}
B & CK_2 \\   
-K_2^{-1}C & K_2^{-1}BK_2 \\
\end{matrix}
\right).
\end{equation*}
In addition, $\det A$ is a sum of two squares over $\FF$.  Specifically, if $i=\sqrt{-1}$, then over 
$\FF[i]$, we have $\det A=\det(B+iC)\det(B-iC)$. \label{SumofSquares}
}
{\pf
Calculating the conditions for which a matrix
\begin{equation*}
A=
\left(
\begin{matrix}
A_1 & A_2 \\  
A_3 & A_4 \\
\end{matrix}
\right)
\end{equation*}
is pseudo-centrosymmetric with respect to $K$ for matrices $K$ of the form given in Equation \eqref{SimpleK} gives us 
that $K_2^{-1}A_2=-A_3K_2$ and $A_4=K_2^{-1}A_1K_2$.  This means we can write
\begin{equation*}
A=
\left(
\begin{matrix}  
B & CK_2 \\
-K_2^{-1}C & K_2^{-1}BK_2 \\
\end{matrix}  
\right),   
\end{equation*}  
for $B=A_1$ and $C=A_2K_2^{-1}$.  With this rewriting, $\det A$ is simple to compute, after adjoining 
$i=\sqrt{-1}$ to $\FF$:
\begin{eqnarray} \nonumber
\det A & = & 
\det   
\left(
\begin{matrix}
B & CK_2 \\
-K_2^{-1}C & K_2^{-1}BK_2 \\
\end{matrix}
\right) \\ \nonumber
& = & 
\det\left(    
\begin{matrix}
I & 0 \\
+iK_2^{-1} & I \\
\end{matrix}  
\right) \det
\left(
\begin{matrix}
B & CK_2 \\
-K_2^{-1}C & K_2^{-1}BK_2 \\
\end{matrix}
\right) \det
\left(
\begin{matrix}
I & 0 \\
-iK_2^{-1} & I \\
\end{matrix}
\right)\\ \nonumber 
& = & \det
\left(
\begin{matrix}
B-iC & CK_2 \\
 0 & K_2^{-1}(B+iC)K_2 \\
\end{matrix}
\right) \\ \nonumber
& = &
\det(B+iC)\det(B-iC),
\end{eqnarray}
as desired.

This is a product of $\det(B+iC)$ and its conjugate, so the determinant of such a matrix $A$ 
is the sum of two squares over $\FF$.
}

In the case of a matrix $A$ that is pseudo-skew-centrosymmetric with respect to $K$, the analogous result states
that $A$ has the form
\begin{equation*}
A=
\left(
\begin{matrix}
B & CK_2 \\ 
K_2^{-1}C & -K_2^{-1}BK_2 \\
\end{matrix}  
\right),   
\end{equation*}
and that $\det A = (-1)^{k}\det(B+iC)\det(B-iC)$, which is also a sum of two squares over $\FF$, possibly up to a sign. The 
sign term appears because the block matrix in the determinant calculation is of the form 
$$\left(
\begin{matrix}
B-iC & CK_2 \\
 0 & -K_2^{-1}(B+iC)K_2 \\
\end{matrix}
\right),$$
with a negative sign appearing from each of the last $k$ rows.  Notice that we can now explicitly evaluate the
determinant of this type of $2k\times 2k$ pseudo- and pseudo-skew-centrosymmetric matrix via a smaller
$k\times k$ determinant.

In Section \ref{sec:general}, we will see that the determinants of all even-order pseudo- and
pseudo-skew-centrosymmetric matrices $A$ can be written as a sum of squares (up to a sign) over the base 
field $\FF$, when $-1$ is not a square in $\FF$.

\subsection{Alternating centrosymmetric matrices}
\label{sec:AltCentro}

We now consider the specific case when the $2k\times 2k$ matrix $K$ is the {\em alternating exchange matrix}---the 
matrix with its cross-diagonal populated with alternating $1$'s and $-1$'s, starting in the upper-right 
corner.  Such a matrix $K$ is anti-involutory.  (Had $K$ been square of odd order, this matrix would have been 
involutory instead of anti-involutory.)

{\defn Let $K$ be the alternating exchange matrix.  We define an $n\times n$ matrix $A$ defined over 
$\FF$ to be {\em 
alternating centrosymmetric} with respect to $K$ if $KA=AK$ and {\em alternating skew-centrosymmetric} with respect
  to $K$ if $KA=-AK$.}

An equivalent classification of $n\times n$ alternating centrosymmetric matrices is that their entries satisfy 
$a_{i,j}=(-1)^{i+j}a_{n+1-i,n+1-j}$.  (See Figure \ref{AltCentroMtx}). An 
$n\times n$ alternating skew-centrosymmetric matrix has entries that satisfy 
$a_{i,j}=(-1)^{i+j+1}a_{n+1-i,n+1-j}$.

\begin{figure}
{\small
\begin{equation*}
\left(
\begin{tabular}{rrrrrr}
$a_1\,$ & $a_2\,$ & $a_3\,$ & $a_4\,$ & $a_5\,$ & $a_6\,$ \\
$a_7\,$ & $a_8\,$ & $a_9\,$ & $a_{10}$ & $a_{11}$ & $a_{12}$ \\
$a_{13}$ & $a_{14}$ & $a_{15}$ & $a_{16}$ & $a_{17}$ & $a_{18}$ \\
$-a_{18}$ & $a_{17}$ & $-a_{16}$ & $a_{15}$ & $-a_{14}$ & $a_{13}$ \\
$a_{12}$ & $-a_{11}$ & $a_{10}$ & $-a_{9}\,$ & $a_{8}\,$ & $-a_{7}\,$ \\
$-a_6\,$ & $a_5\,$ & $-a_4\,$ & $a_3\,$ & $-a_2\,$ & $a_1\,$ 
\end{tabular}
\right)
\end{equation*}
}
\caption{The general form of a $6\times 6$ alternating centrosymmetric matrix}
\label{AltCentroMtx}
\end{figure}

By Theorem \ref{SumofSquares}, we know that the determinant of an alternating centrosymmetric matrix is the 
sum of two squares in $\FF$.  We now present a different version of its determinant with additional 
symmetry built in.

 We first define a set of $k$-member subsets $\tld{I}$ of $[2k]:=\{1,2,\hdots,2k\}$. For any subset $I$ of 
$[k]$, create $\tld{I}$ by taking $I\cup I'$, where $i\in I'$ if $2k+1-i\in [k]\setminus I$.  In this way, 
each $\tld{I}$ has $k$ members.  We will call subsets of $[2k]$ of this form {\em complementary}.  We define 
the sets $S$, $S'$, $T$, and $T'$ of $k$-member subsets of $[2k]$.  If $I$ has $\ell$ elements,
\begin{equation*}
\text{place $\tld{I}$ into set}
\begin{cases} 
S & \text{if $\ell\equiv 0 \mod 4$,} \\
T & \text{if $\ell\equiv 1 \mod 4$,} \\
S' & \text{if $\ell\equiv 2 \mod 4$, or } \\ 
T' & \text{if $\ell\equiv 3 \mod 4$.} \\
\end{cases}
\end{equation*}

Given a $2k\times 2k$ matrix $N$, we define $M(\tld{I})$ to be the $k\times k$ submatrix of $N$ with columns 
restricted to $j\in \tld{I}$, and rows restricted to the first $k$ rows of $N$.

{\thm The formula for the determinant of an alternating centrosymmetric matrix $A$ satisfies 
\begin{equation*}
\det A\! =\!\!
\left[ 
\sum_{\tld{I}\in S} \det(M(\tld{I})) -\!\!
\sum_{\tld{I}\in S'} \det(M(\tld{I})) 
\right]^2
\!\!+
\left[ 
\sum_{\tld{I}\in T} \det(M(\tld{I})) -\!\!
\sum_{\tld{I}\in T'} \det(M(\tld{I})) 
\right]^2\!\!\!.
\label{eq:conjecture}
\end{equation*}
}

{\pf
Let $K_2$ be the upper-right $k\times k$ submatrix of the alternating exchange matrix $K$.
As in the proof of Theorem \ref{SumofSquares}, we consider $\FF[i]$ and write $\det A = 
\det(B+iC)\det(B-iC)$, where $B=A_1$ and $C=A_2K_2^{-1}$.  Calculating $\det(B+iC)$ gives 
$x+iy$ for some $x,y\in\FF$; we will calculate $x$ and $y$ directly.

\begin{equation*}
B+iC=
\left(
\begin{matrix}
        a_{1,1} + ia_{1,2k} & a_{1,2} - ia_{1,2k-1} & \hdots & a_{1,k} + (-1)^{k+1}ia_{1,k+1} \\
        a_{2,1} + ia_{2,2k} & a_{2,2} - ia_{2,2k-1} & \hdots & a_{2,k} + (-1)^{k+1}ia_{2,k+1} \\
              \vdots        &      \vdots           &        &          \vdots                \\
        a_{k,1} + ia_{k,2k} & a_{k,2} - ia_{k,2k-1} & \hdots & a_{k,k} + (-1)^{k+1}ia_{k,k+1} \\
\end{matrix}
\right).
\end{equation*}

Define $b_{j}$ to be the column $b_{j}=(a_{1,j},a_{2,j},\hdots,a_{k,j})^{T}$.  By linearity of determinants, 
$\det (B+iC)$ is the sum of $2^k$ determinants of matrices $M$ with dimensions $k\times k$, where in column 
$j$ we can choose to place either $b_{j}$ or $i(-1)^{j+1}b_{2k+1-j}$.  Given any determinant of this form, we 
can convert it into a form where the indices of the columns are increasing: 
$\tld{I}=i_1<\cdots<i_r<k+1/2<i_{r+1}<\cdots<i_k$. Note that $\tld{I}$ is complementary as defined above. When 
we do this and account for changes of sign by interchanging columns, the matrix $M$ becomes
\begin{equation*}
M=
\left(
\begin{matrix}
        a_{1,i_1} & \hdots &  a_{1,i_r} & i(-1)^{k+1}a_{1,i_{r+1}} & \hdots & i(-1)^{k+1}a_{1,i_k} \\
        a_{2,i_1} & \hdots &  a_{2,i_r} & i(-1)^{k+1}a_{2,i_{r+1}} & \hdots & i(-1)^{k+1}a_{2,i_k} \\
         \vdots   &        &   \vdots   &     \vdots     &        &   \vdots   \\
        a_{k,i_1} & \hdots &  a_{k,i_r} & i(-1)^{k+1}a_{k,i_{r+1}} & \hdots & i(-1)^{k+1}a_{k,i_k} \\
\end{matrix}
\right).
\end{equation*}
The determinant of this matrix is $(i(-1)^{k+1})^{k-r}$ times $\det(M(\tld{I}))$.  In particular, matrices such
that $(k-r)\equiv 0 \mod 2$ contribute to $x$, while matrices such that $(k-r)\equiv 1
\mod 2$ contribute to $y$.  In addition, the value of $(k-r)$ mod 4 determines the
sign of the contribution to the sum.  

This establishes the theorem.  
}

A similar approach gives an analogous statement for alternating skew-centro\-sym\-metric matrices---we have instead,
\begin{equation*}
(-1)^k\det A\! =\!\!
\left[ 
\sum_{\tld{I}\in S} \det(M(\tld{I})) -\!\!
\sum_{\tld{I}\in S'} \det(M(\tld{I})) 
\right]^2
\!\!+
\left[ 
\sum_{\tld{I}\in T} \det(M(\tld{I})) -\!\!
\sum_{\tld{I}\in T'} \det(M(\tld{I})) 
\right]^2\!\!\!.
\label{eq:conjecture2}
\end{equation*}
Just as in the analogous statement of Theorem \ref{SumofSquares}, a $(-1)^k$ term appears.

\subsection{Pseudo-centrosymmetric matrices in general}
\label{sec:general}

An anti-involutory matrix $K$ behaves like $\sqrt{-1}$ when multiplying a vector.  When $-1$ is not a square in
$\FF$, such a $K$ allows us to create a very nice basis of $\FF^{2k}$, similar to 
the construction of an almost complex structure on $\RR^{2k}$.  For more 
background on almost complex structures, see Section 5.2 of Goldberg's Curvature and Homology 
\cite{Goldberg}.  The following lemma will help us prove a general result about even-order pseudo- and 
pseudo-skew-centrosymmetric matrices.

{\lem If $K$ is a $2k \times 2k$ anti-involutory matrix with entries in $\FF$, a field not of characteristic 2 in
  which $-1$ is not a square, then $\FF^{2k}$ has a basis of the form $\{v_1,\hdots,v_n,Kv_1,\hdots,Kv_n\}$ for
  vectors $v_i\in \FF^{2k}$. \label{lem:basis}}
{\pf We will proceed by induction of a set of independent vectors 
$S_l=\{v_1,\hdots,v_l,Kv_1,\hdots,Kv_l\}$.  Given any nonzero vector $v_1\in \FF^{2k}$, we show that 
$S_1=\{v_1,Kv_1\}$ is linearly independent.  Suppose to the contrary that $cv_1=Kv_1$ for some $c\in\FF$.  
then $c^2v_1=K(cv_1)=K^2v_1=-v_1$, so $(c^2+1)v_1=0$, which is impossible since $\sqrt{-1}\notin \FF$. This 
proves the base case.  

Now suppose that $S_{l}$ is a set of $2l$ independent vectors of the form above, with $l<k$.  Since 
$2l<2k,$ there exists a vector $v_{l+1}\in \FF^{2k}$ not in the span of $S_l$.  Now we wish to show that 
$S_{l+1}=S_l\cup\{v_{l+1},Kv_{l+1}\}$ is linearly independent.  Consider $x=v_{l+1}-\proj_{S_l}(v_{l+1})$ and 
$y=Kv_{l+1}-\proj_{S_l}(Kv_{l+1})$, the components of $v_{l+1}$ and $Kv_{l+1}$ orthogonal to $S_l$.  Notice that 
$Kx=y$, so by a similar argument to above, $x$ and $y$ are linearly independent.  Therefore $S_{l+1}$ is 
linearly independent, and the lemma follows by induction. }

This lemma is the starting-off point for the following theorem; we rely on the simpler result from 
Theorem \ref{SumofSquares}.

{\thm
 Let $K$ be a $2k\times 2k$ anti-involutory matrix with entries in $\FF$, a field not of characteristic 2 in which
 $-1$ is not a square.  Let $A$ be a matrix (also with entries from $\FF$) which satisfies $KA=AK$.
  Then the determinant of $A$ is a sum of two squares of elements from $\FF$.
\label{thm:SumofSquaresGen}
}
{\pf
Choose $k$ vectors $\{{\bf v}_i\}_{i=1}^k \in \FF^{2k}$ such that the set $\{{\bf v}_1, \hdots, {\bf 
v}_k, K{\bf v}_1, \hdots, K{\bf v}_k\}$ forms a basis for $\FF^{2k}$, as guaranteed by Lemma \ref{lem:basis}.
Define the $2k\times 2k$ matrix $V=\left[{\bf v}_1|\cdots|{\bf v}_k|K{\bf v}_1|\cdots|K{\bf v}_k\right]$.  
Since the column vectors form a basis for $\FF^{2k}$, the matrix is invertible.  Notice that the product 
$KV=\left[K{\bf v}_1|\cdots|K{\bf v}_k|-{\bf v}_1|\cdots|-{\bf v}_k\right]=VK'$, where 
$K'=\left[\begin{matrix} 0 & -I_k \\ I_k & 0\end{matrix}\right]$.  Therefore, $K=VK'V^{-1}$.

Now, if $KA=AK$, then $VK'V^{-1}A=AVK'V^{-1}$, so $K'(V^{-1}AV)=(V^{-1}AV)K'$.  We see that $K'$ is of the form in
Equation \eqref{SimpleK}, so $\det(A)=\det(V^{-1}AV)$ is a sum of two squares of elements from $\FF$, by 
Theorem \ref{SumofSquares}.
}

When $A$ is a pseudo-skew-centrosymmetric matrix, the argument is the same up until the last sentence, where we
recognize that if $KA=-AK$, $(V^{-1}AV)$ is pseudo-skew-centrosymmetric with respect to $K'$, and the remark after
Theorem \ref{SumofSquares} implies that the determinant of $A$ is a sum of two squares, possibly up to a sign.

{\cor Let $K$ be an anti-involutory matrix over $\QQ$ and let $A$ be a matrix with integral entries satisfying
 $AK=KA$ (or $AK=-KA$).  Then $\det A$ is a sum of two integral squares (possibly up to a
 sign). \label{cor:Integral}} {\pf It is clear that $\det A$ is an integer.  By Theorem \ref{thm:SumofSquaresGen}
 with $\FF=\QQ$, $\det A$ is the sum of two rational squares.  (Possibly up to a sign if $AK=-KA$.)  An integer is
 a sum of two rational squares only if it is a sum of two integral squares. This follows from the characterization
 that an integer $n$ can be written as a sum of two squares if and only if its prime factorization contains only
 even powers of primes of the form $4k+3$ (See \cite{Davenport}, p.116).  If we can write $(p_1/q_1)^2+(p_2/q_2)^2=n$ for integers $p_1$,
 $q_1$, $p_2$, and $q_2$, then $p_1^2+p_2^2=q_1^2q_2^2n$; hence $q_1^2q_2^2n$ is a sum of two integral squares.
 This implies that $q_1^2q_2^2n$ (and $n$ itself) only contain primes of the form $4k+3$ to even powers.  }


\section{Applications to Enumerative Combinatorics}
\label{sec:Enum}

The matrix-theoretical results from Sections \ref{sec:PseudoCentro} and \ref{sec:AltCentro} allow us to reprove
some results from matching theory.

\subsection{Revisiting Jockusch's theorem}
\label{sec:JockuschThm}

A certain symmetry property of a bipartite graph allows us to say something about its number of perfect matchings.

{\defn A {\em 2-even symmetric graph} $G$ is a connected planar bipartite graph such that a 180 degree 
rotation $R_2$ about the origin maps $G$ to itself and the length of a path between $v$ and $R_2(v)$ is even.}

\begin{figure}
\begin{center}
\epsfig{figure=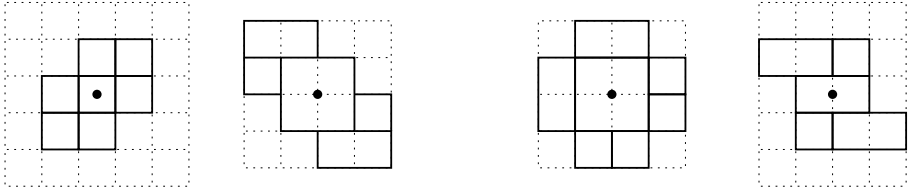}
  \caption{Two graphs that are 2-even symmetric, and two graphs that are not.}
  \label{SymExample}
\end{center}
\end{figure}

In particular, if the graph is embedded in a square grid, if you rotate the graph by 180 degrees, you get the 
same graph back, AND the center of rotation is not the center of an edge.  (See Figure \ref{SymExample}.) 
In \cite{Jockusch}, William Jockusch proves that if a graph is 2-even symmetric then the number of perfect 
matchings of the graph is a sum of squares.  Jockusch's result produces a weighted labeling function $u$ of 
the quotient graph $G_2$ involving complex numbers.  Once one counts the number of weighted matchings 
associated to $u$, denoted $M_u(G_2)$, the number of perfect matchings of $G$ is 
$M_u(G_2)\overline{M_u(G_2)}$, 
resulting in a sum of two squares.

We can use Theorem \ref{SumofSquares} to reprove Jockusch's theorem for a subset of 2-even symmetric graphs 
that relies only on the structure of the Kasteleyn-Percus matrix of the region, to be discussed shortly 
hereafter.  This subset of graphs occurs when the 2-even symmetric graph can be embedded in the square lattice 
with the center of rotation in the center of a square.  When we restrict our graphs to this type, we can 
represent the square lattice in the standard $x$-$y$ coordinates by placing vertices at 
$(2k+1,2\ell+1)$ for $k,\ell\in \ZZ$, so that the center of rotation $(0,0)$ is the centroid of some square in 
the lattice.  We color the vertices white if $k+\ell$ is even and black if $k+\ell$ is odd.  In Section 
\ref{sec:CtrEx}, we present an example that shows that a restriction of this type is necessary.

\begin{figure}
\begin{center}
\epsfig{figure=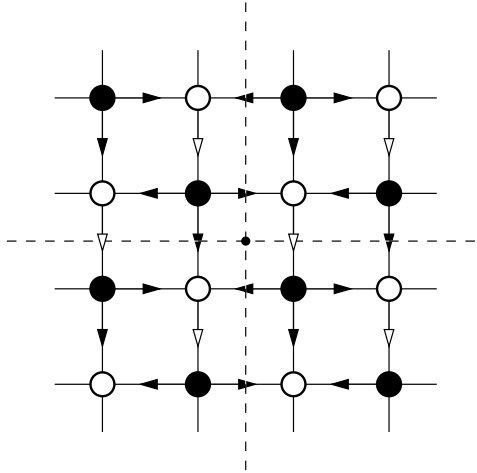}
  \caption{The canonical orientation of edges on the square lattice}
  \label{Orientations}
\end{center}
\end{figure}

Recall that a bipartite graph $G=(V,E)$ is a graph where the vertex set is partitioned into two subsets, the
``black'' and ``white'' vertices, where there are no edges connecting vertices of the same color.  Notice that if
there is to be a perfect matching (or pairing of all vertices using graph edges), there must be the same number of
black vertices as white vertices.  A well-known method of counting perfect matchings of a planar bipartite graph
relies on taking the determinant of a Kasteleyn-Percus matrix.  (See \cite{Kasteleyn, Percus} for background.)  The
definition we will use in this article is that the {\em Kasteleyn-Percus matrix} $A$ of a bipartite graph $G=(V,E)$
has $|V|/2$ rows representing the white vertices and $|V|/2$ columns representing the black vertices.  The non-zero
entries $a_{ij}$ of $A$ are exactly those that have an edge between white vertex $w_i$ and black vertex $b_j$.
These entries are all $+1$ or $-1$ depending on the position of the edges they represent in the graph---the
restriction is that elementary cycles have a net $-1$ product.  In the case of the square lattice above, we can
satisfy this condition easily by giving matrix entries the value $-1$ if they correspond to edges that are of the
form $e=(v_1\,v_2)$ with $v_1=(2k-1,2\ell+1)$ and $v_2=(2k+1,2\ell+1)$ and such that $v_1$ is black.  This is most
easily understood by giving orientations to the edges of the lattice as in Figure \ref{Orientations}, and assigning
an edge the value $+1$ if the edge goes from black to white and the value $-1$ if the edge goes from white to
black.  The absolute value of the determinant of this matrix counts the number of perfect matchings of $G$.

With this definition of the Kasteleyn-Percus matrix, we can formulate the following theorem.

{\thm
Let $G$ be a $2$-even symmetric graph embedded in the square lattice with the rotation axis in 
the center of one square. Under a suitable ordering of the vertices, the Kasteleyn-Percus matrix $A$ of $G$ is 
alternating centrosymmetric. \label{JockuschExt} 
}

{\pf
We label the black and white vertices to determine the positions of the $+1$ and $-1$ 
entries in 
$A$.  After an initial labeling, we interchange rows and columns as necessary to manipulate the matrix into 
being alternating centrosymmetric, as follows.

Embedded in this lattice, half the vertices of $G$ lie above the horizontal line through the origin.  
Coloring the vertices of $G$ the color they inherit from the lattice coloring above, $R_2$ takes vertices to 
counterparts of the same color so for some $m$, we have $m$ vertices of each color in the upper half of the 
graph and $4m$ vertices in all.  Label all white vertices $v$ in the upper half of graph with values $1$ to $m$, 
and do the same for black vertices $w$.  For each vertex $x$ with value $i$, label $R_2(x)$ with value $2m+1-i$. 
(See Figure \ref{Labeling} for an example.)

\begin{figure}
\begin{center}
\begin{tabular}{ccccccc}
  \epsfig{figure=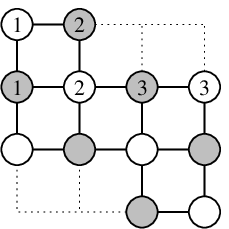} &\!\!&
  \epsfig{figure=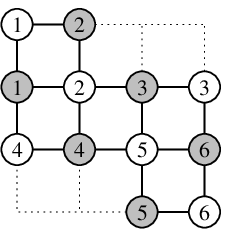} &\!\!&
  \epsfig{figure=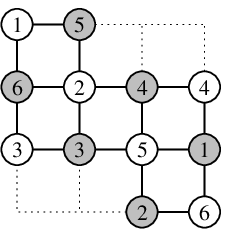} &\!\!\!&
\begin{tabular}{c}\vspace{.8in}{\small $\left[\begin{array}{cccccc}
0 & 0 & 0 & 1 & 1 & \!\!-1\!\! \\
0 & 0 & 0 & 0 & 1 & 1 \\
0 & 1 & 1 & 0 & 1 & 0 \\
0 & 1 & 0 & 1 & \!\!-1\!\! & 0 \\
1 & \!\!-1\!\! & 0 & 0 & 0 & 0 \\
1 & 1 & \!\!-1\!\! & 0 & 0 & 0 \\
\end{array}\right]$}
\end{tabular}
\end{tabular}\vspace{-.8in}
  \caption{Example of the labeling procedure in the proof of Theorem \ref{JockuschExt}.}
  \label{Labeling}
\end{center}
\end{figure}

From this initial numbering of vertices, we wish to modify some labels so that the labels in each row of the 
square lattice are of the same parity.  Note that each vertex and its counterpart have opposite parity.  
Start with the top row.  For each vertex that is labeled with an even number, switch its label with its 
counterpart.  In this way, all elements in the top row will have an odd value.  For the second row, exchange a 
vertex's label with its counterpart's label if it has an odd value.  Continue in this fashion until all odd rows 
have vertices with odd labels and all even rows have vertices with even values.  After determining all rows 
above the horizontal line through the origin, the rest of the rows come for free.

By this construction, for any horizontal edge $(v_i,w_j)$ in $G'$, we know that $i$ and $j$ are of the same 
parity.  Similarly, we know that any vertical edge $(v_i,w_j)$ has $i$ and $j$ of opposite parity. In 
addition, this implies that the rotation of a vertical edge by $R_2$ results in the opposite sign appearing in 
the Kasteleyn-Percus matrix.

Since $(v_i,w_j)$ is a horizontal edge if and only if $(v_{2m+1-i},w_{2m+1-j})$ is a horizontal edge, a $+1$ 
appears in position $a_{(i,j)}$ for $i+j$ even if and only if $a_{(2m+1-i,2m+1-j)}$ is $+1$.  Similarly, if 
$(v_i,w_j)$ is a vertical edge, then so is $(v_{2m+1-i},w_{2m+1-j})$, and their entries in $A$ are opposite.  
This occurs exactly where $i+j$ is odd.  All other entries are zero, so for those $a_{i,j}$ we have 
$a_{(i,j)}=\pm a_{(2m+1-i,2m+1-j)}$.

These conditions imply that entry $a_{(i,j)}$ equals $(-1)^{i+j}a_{(2m+1-i,2m+1-j)}$, which implies $A$ is 
alternating centrosymmetric, as desired.  
}

With this theorem and Corollary \ref{cor:Integral}, we have the following corollary.

{\cor The number of perfect matchings of a 2-even-symmetric graph embedded in the square lattice with
  the center of rotation in the center of one unit square is a sum of two integral squares. \label{SOS}
}

{\rem This theorem is weaker than Jockusch's original theorem as it applies to fewer regions. 
However, it allow us to prove a sum-of-squares result using only matrix-based methods.
}

\subsection{Applications to generalized Aztec pillows}
\label{sec:Pillows}

\begin{figure}
\begin{center}
\begin{tabular}{ccc}
  \epsfig{figure=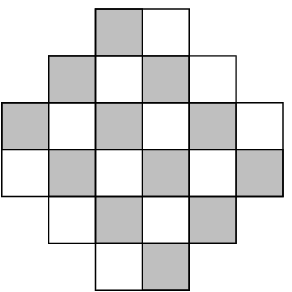} &
  \epsfig{figure=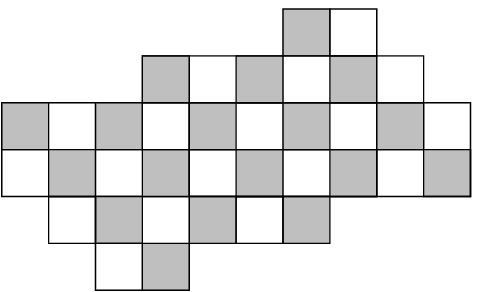} &
  \epsfig{figure=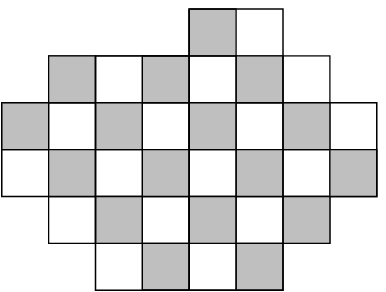}
\end{tabular}
  \caption{Examples of an Aztec diamond, an Aztec pillow, and a generalized Aztec pillow}
  \label{AztecEx}
\end{center}
\end{figure}

\begin{figure}
\begin{center}
\begin{tabular}{ccc}
  \epsfig{figure=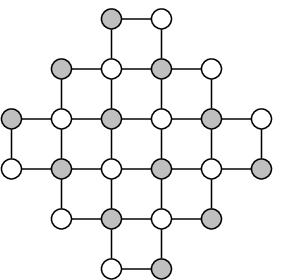} &
  \epsfig{figure=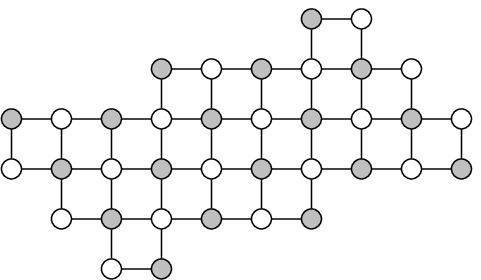} &
  \epsfig{figure=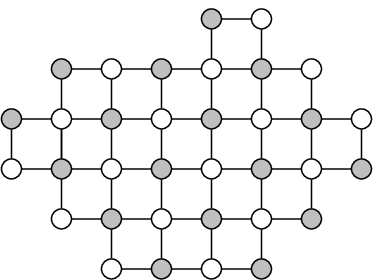}
\end{tabular}
  \caption{The dual graphs to the regions in Figure \ref{AztecEx}.}
  \label{AztecEx2}
\end{center}
\end{figure}

An application of Corollary \ref{SOS} has to do with domino tilings of rotationally-symmetric regions made up of 
unit squares.  A {\em domino tiling} of a region is a complete covering of the region with non-overlapping $2\times 1$ 
and $1\times 2$ rectangles (or dominoes).  We can associate to any region its {\em dual graph}, a graph with vertices
representing the unit squares and edges between vertices representing adjacent squares.  In effect, counting the 
number of domino tilings of a given region is the same as counting the number of perfect matchings of the 
corresponding dual graph.  Many results in this vein are presented in James Propp's survey article \cite{Propp}.  
Some regions fit the framework from the previous section particularly well; examples are pictured in Figure 
\ref{AztecEx} and their corresponding dual graphs are pictured in Figure \ref{AztecEx2}.  

An {\em Aztec diamond} is the union of the $2n(n+1)$ unit squares with integral vertices $(x,y)$ such that
$|x|+|y|\leq n+1$.  {\em Aztec pillows} were introduced in \cite{Propp} and explored more in depth in
\cite{HamThm}.  They are also rotationally-symmetric regions composed of unit squares with their ``steps'' along
the northwest and southeast diagonals having height one and length three.  Lastly, we include {\em generalized
Aztec pillows}, where all steps off the central band of squares are of height one and odd length.  Notice that this
implies that Aztec diamonds and regular Aztec pillows are also generalized Aztec pillows.  While generalized Aztec
pillows need not be rotationally symmetric, when we restrict to those that are rotationally symmetric, we now have
many regions whose dual graphs are 2-even-symmetric graphs, so we have the following corollary.

{\cor
The number of domino tilings of any rotationally-symmetric generalized Aztec pillow, which 
includes all 
Aztec diamonds and Aztec pillows, is a sum of two integral squares.  \label{APsumofSquares}
}

\subsection{Problems extending Theorem \ref{JockuschExt}}
\label{sec:CtrEx}

\begin{figure}
\begin{center}
\epsfig{figure=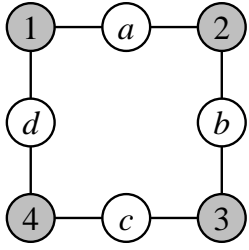}
  \caption{A 2-even symmetric graph whose Kasteleyn-Percus matrix is not alternating centrosymmetric.}
  \label{LabelingCtrEx}
\end{center}
\end{figure}

Unfortunately, the methods from Section \ref{sec:JockuschThm} do not allow us to completely reprove Jockusch's
theorem, as we highlight with the following example.  Consider the graph in Figure \ref{LabelingCtrEx}.  This graph
is clearly bipartite, can be embedded in the square grid, and is 2-even symmetric.  However, the following theorem
holds.

\medskip
{\thm
No reordering of the vertices transforms the Kasteleyn-Percus matrix of the graph in Figure \ref{LabelingCtrEx} 
into an alternating centrosymmetric matrix.
}

{\pf
Notice that our relabeling trick from Section \ref{sec:JockuschThm} will not work here, as the $x$-axis would 
contain both even- and odd-labeled vertices.  In addition, any rearrangement of the vertices could produce 
neither an alternating centrosymmetric matrix nor an alternating skew-centrosymmetric matrix.  We can see this 
as follows.  Vertex $a$ is adjacent to vertices 1 and 2, while vertex $c$ is adjacent to vertices 3 and 4, all 
by $+1$-weighted edges through our scheme from Figure \ref{Orientations}.  This implies that whatever rows of 
our rearranged matrix $A$ the vertices $a$ and $c$ are in, the entries used in the rows will be complementary 
(such as of type (I) 1010 and 0101, type (II) 1001 and 0110, or type (III) 1100 and 0011).  
This is true also with $b$ and $d$, with the addition of some signs.  Since vertices $b$ and $d$ share one 
vertex each with $a$ and $c$, they can not have complementary pairs of the same type (I, II, or III).  In 
order for the non-zero entries of $A$ to match up correctly, $a$'s and $c$'s rows must be the center two rows 
or the first and last rows.  Therefore only types (I) and (III) are valid types, and there must be one of each 
for the two pairs.  Unfortunately, this can not possibly work when considering the sign conventions necessary 
for a matrix to be alternating centrosymmetric or alternating skew-centrosymmetric. 
}

{\rem As alluded to in the introduction of Section \ref{sec:JockuschThm}, there are multiple definitions of 
a Kasteleyn-Percus matrix.  Perhaps it is possible to reprove Jockusch's theorem using another matrix 
interpretation.
}


{\rem Although we may not be able to reprove Jockusch's theorem using matrix methods for all 2-even 
symmetric graphs, perhaps the condition that the graph be embedded in the square lattice can be relaxed.
}

\section{Acknowledgments}

I would like to thank Mark Yasuda for his corrections and for helping strengthen the focus of this paper.  I also thank
Henry Cohn for his ideas, critiques, and support.  I thank Paul Loya for his help with almost complex structures.
A distant draft of this article appeared in my doctoral dissertation.  

\nocite{*}   
\bibliographystyle{elsart-num}
\bibliography{AltCentroSymm}

\end{document}